\documentclass{article}
\usepackage{graphicx}

\begin{document}

\section*{A quantitative program for Hadwiger's covering conjecture and
Borsuk's partition conjecture\footnote{This work is supported by the
National Science Foundation of China, the Chang Jiang Scholars
Program and LMAM at Peking University.}}

\noindent {\bf Chuanming ZONG}

\bigskip \noindent School of Mathematical Sciences, Peking
University, Beijing 100871, China

\noindent (Email: cmzong@math.pku.edu.cn)

\bigskip \noindent Dedicated to Professor Yuan Wang on the occasion
of his 80th birthday

\bigskip
\noindent {\bf Abstract:}  In 1957, Hadwiger made a conjecture that
every $n$-dimensional convex body can be covered by $2^n$ translates
of its interior. Up to now, this conjecture is still open for all
$n\ge 3$. In 1933, Borsuk made a conjecture that every
$n$-dimensional bounded set can be divided into $n+1$ subsets of
smaller diameters. Up to now, this conjecture is open for $4\le n\le
297$. In this article we encode the two conjectures into continuous
functions defined on the spaces of convex bodies, propose a
four-step program to approach them, and obtain some partial results.

\medskip \noindent {\bf Keywords:} convex body, Hadwiger's
conjecture, Banach-Mazur metric, $\beta$-net, Borsuk's conjecture

\medskip \noindent  {\bf MSC (2000):}  52C17; 11H31

\section*{1. Introduction}

\medskip
In $n$-dimensional {\it Euclidean space} $E^n$, let $K$ be a {\it
convex body} with boundary $\partial (K)$, interior $int (K)$ and
volume $v(K)$, and let $c(K)$ denote the smallest number of
translates of $int (K)$ that their union can cover $K$. In 1955,
Levi \cite{levi54} studied $c(K)$ for the two-dimensional convex
domains and proved that
$$c(K)=\left\{
\begin{array}{ll}
4,& \mbox{if $K$ is a parallelogram,}\\
3,& \mbox{otherwise.}
\end{array}\right.$$
Let $P$ denote an $n$-dimensional {\it parallelopiped}. Clearly, any
translate of $int (P)$ can not cover two vertices of $P$. Therefore,
it can be deduced that
$$c(P)=2^n.$$
Let $B$ denote an $n$-dimensional ball. It is easy to see that
$$c(B)=n+1.$$
In fact, this is true for all $n$-dimensional convex bodies with
smooth boundaries.

Based on these results and some other observations, in 1957 Hadwiger
\cite{hadw57} made the following conjecture: {\it For every
$n$-dimensional convex body $K$ we have
$$c(K)\le 2^n,\eqno (1)$$
where the equality holds if and only if $K$ is a parallelopiped.}

This conjecture has been studied by many authors. In the course,
many partial results have been achieved and several connections with
other important problems such as the {\it illumination problem} and
the {\it separation problem} have been discovered (see Bezdek
\cite{bezd06}, Boltyanski, Martini and Soltan \cite{bolt97}, Brass,
Moser and Pach \cite{bras05} and Zong \cite{zong08} for general
references). For example, Lassak \cite{lass84} proved this
conjecture for the three-dimensional centrally symmetric case,
Rogers and Zong \cite{roge97} obtained
$$c(K)\le {2n\choose n}(n\log n+n\log \log n+5n)$$
for general $n$-dimensional convex bodies, and
$$c(K)\le 2^n (n\log n+n\log \log n+5n)$$
for centrally symmetric ones. Nevertheless, we are still far away
from the solution of the conjecture, even the three-dimensional
case.

Let $m$ be a positive integer and let $\gamma_m(K)$ be the smallest
positive number $r$ such that $K$ can be covered by $m$ translates
of $rK$. Clearly, we have
$$\gamma_m(K)=1$$
for all $m\le n$, and
$$\gamma_m(K)\ge \gamma_{m+1}(K)$$
for all positive integers $m$ and all $n$-dimensional convex bodies
$K$.

Let $\mathcal{T}^n$ denote the set of all non-singular linear
transformations in $E^n$ and let $\mathcal{K}^n$ denote the space of
all $n$-dimensional convex bodies with the {\it Banach-Mazur metric}
defined by
$$\| K_1, K_2\| =\log\ \min \left\{ r: \ K_1\subseteq T(K_2)\subseteq
rK_1+{\bf x}; \ {\bf x}\in E^n; \ T\in \mathcal{T}^n\right\}.$$ By
John's theorem (see Section 3) and Blaschke's selection theorem (see
\cite{grub07}) it follows that $\mathcal{K}^n$ is bounded, connected
and compact. On the other hand, for any given positive integer $m$,
it can be shown that $\gamma_m(K)$ as a function of $K$ defined on
$\mathcal{K}^n$ is continuous (see Theorem A in Section 3). In
addition, we have
$$\gamma_m(K_1)=\gamma_m(K_2)$$
whenever $\| K_1, K_2\| =0$. Then we define $$\Gamma (n,m)
=\max_{K\in \mathcal{K}^n}\{ \gamma_m(K)\}$$ and
$$\gamma (n,m)
=\min_{K\in \mathcal{K}^n}\{ \gamma_m(K)\}.$$

Note that $\gamma_m(K)$ is continuous on $\mathcal{K}^n$ and
$\mathcal{K}^n$ is compact. It is easy to see that (1) holds for all
$n$-dimensional convex bodies $K$ if and only if
$$\gamma_{2^n}(K)<1$$
holds for all $K\in \mathcal{K}^n$. Therefore, it is equivalent to
$$\Gamma (n,2^n)<1.$$
Thus, Hadwiger's covering conjecture can be encoded into the
functions $\gamma_m(K)$ defined on the space $\mathcal{K}^n$.

In this article, based on some nice properties of $\gamma_m(K)$ and
$\mathcal{K}^n$, we suggest a four-step program (see Section 3) to
approach Hadwiger's conjecture as well as Borsuk's conjecture. In
addition, we study the values of $\gamma_m(K)$ for some particular
$m$ and $K$. Among other things, the following results are proved:

\medskip\noindent
{\bf Theorem 1.} {\it Let $K$ be a bounded three-dimensional convex
cone $($the convex hull of a convex domain and a point which is not
in the plane of the domain$)$, then we have}
$$\gamma_8(K)\le \mbox{${2\over 3}$}.$$

\medskip\noindent
{\bf Theorem 2.} {\it Let $K_p$ be the unit ball of the
three-dimensional $\ell_p$ norm, in other words, $$K_p=\left\{ (x,
y, z):\ |x|^p+|y|^p+|z|^p\le 1\right\}.$$ For all $p$ satisfying
$1\le p\le +\infty$ we have}
$$\gamma_8(K_p)\le \sqrt{\mbox{$2\over 3$}}.$$

\section*{2. The two-dimensional case, a brief review}

\medskip
The values of $\gamma (2, m)$ and $\Gamma (2, m)$ have been studied
by several authors. Clearly, we have $$\gamma (2,2)=\Gamma
(2,2)=\Gamma (2,3)=1$$ and, by considering the area measures,
$$\gamma (2, m)\ge {1\over {\sqrt{m}}}.$$

However, for the nontrivial cases, it is not easy to determine the
exact values of $\gamma (2, m)$ and $\Gamma (2, m)$. We list the
known results in Table 1 and Table 2.
\begin{table}[ht]
\renewcommand\arraystretch{1.5}
\noindent\[
\begin{array}{|c|c|c|c|}
\hline
m&3&4&5\\
\hline
\gamma (2,m)&{2\over 3}&{1\over 2}&{1\over 2}\\
\hline
{\rm Authors}&{\rm J.F.\ Belousov\ \cite{belo77}}&{\rm S.\ Krotoszynski\ \cite{krot87}} &
{\rm S.\ Krotoszynski\ \cite{krot87}}\\
\hline
\end{array}
\]
\centerline{Table 1}
\end{table}

\begin{table}[ht]
\renewcommand\arraystretch{1.5}
\noindent\[
\begin{array}{|c|c|c|c|c|c|c|}
\hline
m&3&4&5&6&7&8\\
\hline
\Gamma (2,m)&1&{\sqrt{2}\over 2}&??&??&{1\over 2}&{1\over 2}\\
\hline {\rm Authors}&&{\rm M.\ Lassak\ \cite{lass86}}&{\rm ??} &{\rm
??} &{\rm F.W.\ Levi\ \cite{levi54}}
&{\rm F.W.\ Levi\ \cite{levi54}}\\
\hline
\end{array}
\]
\centerline{Table 2}
\end{table}

\medskip\noindent {\bf Remark 1.} By $\Gamma (2,2^2)=\sqrt{2}/2$ it
follows that every two-dimensional convex domain $K$ can be covered
by four translates of ${\sqrt{2}\over 2}K$. As shown in Table 2, the
values of $\Gamma (2,5)$ and $\Gamma (2,6)$ have not been determined
yet.

\section*{3.  A four-step program for Hadwiger's conjecture and Borsuk's conjecture}

\medskip
First, let us introduce a basic result about $\gamma_m(K)$.

\medskip\noindent
{\bf Theorem A.} {\it For any pair of positive integer $m$ and $n$,
the function $\gamma_m(K)$ is continuous on $\mathcal{K}^n$.}

\medskip\noindent
{\bf Proof.} Let $\lambda $ and $r$ be positive numbers, $r\le 1$,
let $K$ be an $n$-dimensional convex body, and let ${\bf x}_i$ be
$m$ suitable points. Assume that
$$\gamma_m(K)=r$$
and
$$K\subseteq \bigcup_{i=1}^m\left( rK+{\bf x}_i\right).$$

For each point ${\bf x}\in K$, there is a corresponding point ${\bf
y}\in K$ and an index $i$ satisfying $${\bf x}=r{\bf y}+{\bf x}_i.$$
Thus we can deduce
$$\lambda {\bf x}=\lambda r {\bf y}+\lambda {\bf x}_i$$
and
$$\lambda K\subseteq \bigcup_{i=1}^m\left( \lambda r K+\lambda {\bf
x}_i\right).\eqno (2)$$

Let $\epsilon $ be a small positive number and let $K'$ be any
$n$-dimensional convex body satisfying
$$\| K, K'\| \le \log (1+\epsilon ).$$
In other words, without loss of generality, we may take
$$K\subseteq K'\subseteq (1+\epsilon )K.$$
Then, by (2) we have
\begin{eqnarray*}
K' &\subseteq & (1+\epsilon )K\subseteq \bigcup_{i=1}^m\left(
(1+\epsilon )r K+(1+\epsilon ){\bf x}_i\right)\\
& \subseteq &\bigcup_{i=1}^m\left( (1+\epsilon )r K'+(1+\epsilon
){\bf x}_i\right), \end{eqnarray*} which implies that
$$\gamma_m(K')\le r+\epsilon r\le \gamma_m(K)+\epsilon .$$
Similarly, we can get
$$\gamma_m(K)\le \gamma_m(K')+\epsilon $$
and therefore
$$| \gamma_m(K')-\gamma_m(K)|\le \epsilon ,$$
which means that the function $\gamma_m(K)$ is continuous at $K$.
The theorem is proved. \hfill{$\diamondsuit $}

\medskip\noindent
{\bf Remark 2.} In fact, it follows from the proof that
$\gamma_m(K)$ is uniformly continuous on $\mathcal{K}^n$.

\medskip
Let $B^n$ denote the $n$-dimensional unit ball centered at
the origin. In 1948, F. John \cite{john48} proved the following
result:

\medskip
\noindent {\bf John's theorem.} {\it For each $n$-dimensional convex
body $K$ there is an non-singular linear transformation $T\in
\mathcal{T}^n$ such that}
$$B^n\subseteq T(K)\subseteq nB^n.$$

\medskip
Let $\overline{\mathcal{K}^n}$ denote the set of all convex bodies
$K$ satisfying
$$B^n\subseteq K\subseteq nB^n.\eqno (3)$$
By John's theorem we have
$$\Gamma (n,m)=\max_{K\in \overline{\mathcal{K}^n}}\gamma_m(K).$$

\noindent {\bf Definition 1.} Let $\beta $ be a positive number, and
let $K_1,$ $K_2$, $\cdots$, $K_{l(n,\beta )}$ be $l (n,\beta )$
convex bodies in $\mathcal{K}^n$, where $l (n,\beta )$ is an integer
depending on $n$ and $\beta $. If for any $K\in \mathcal{K}^n$ there
is a corresponding $K_i$ satisfying
$$\| K, K_i\| \le \beta ,$$
then we call $\mathcal{N}=\{ K_1, K_2, \cdots , K_{l(n,\beta )}\}$ a
$\beta$-net in $\mathcal{K}^n$.

\medskip
\noindent {\bf Remark 3.} Defining
$$\mathcal{B}(K_i,\beta )=\left\{ K\in \mathcal{K}^n: \ \| K,
K_i\|\le \beta \right\},$$ it is easy to see that $\mathcal{N}=\{
K_1, K_2, \cdots , K_{l(n,\beta )}\}$ is a $\beta $-net in
$\mathcal{K}^n$ if and only if
$$\bigcup_{i=1}^{l(n,\beta )}\mathcal{B}(K_i, \beta )=\mathcal{K}^n.$$

The existence of the $\beta $-nets in $\mathcal{K}^n$ is guaranteed
by the following lemma which seems important in the study of
$\mathcal{K}^n$.

\medskip
\noindent {\bf Fundamental Lemma.} {\it For each pair $\{ n, \beta
\}$ of positive integer $n$ and positive number $\beta $ there is a
corresponding integer $l (n,\beta )$ such that $\mathcal{K}^n$ has a
$\beta $-net of $l (n,\beta )$ elements.}

\medskip\noindent
{\bf Proof.} Let $\theta $ be a small positive number and assume
that
$$X=\left\{ {\bf x}_1, {\bf x}_2, \cdots , {\bf x}_{c(n,\theta
)}\right\}$$ is a subset of $\partial (nB^n)$ such that the
$c(n,\theta )$ caps $(\theta B^n+{\bf x}_i)\cap \partial (nB^n)$
form a covering on $\partial (nB^n)$. For convenience, we use
$\theta '$ to denote the spherical radii of the caps.

Let $m$ be a large integer and define
$$X_{i,m}=\left\{ \mbox{$1\over n$}{\bf x}_i+\mbox{${{j(n-1)}\over
{mn}}$}{\bf x}_i:\ j=0, 1, \cdots , m\right\}$$ and
$$\mathcal{P}=\left\{ conv\{ {\bf y}_1, {\bf y}_2, \cdots , {\bf
y}_{c(n,\theta )}\}:\ {\bf y}_i\in X_{i, m}\right\}.$$ Note that
${1\over n}{\bf x}_i\in \partial (B^n)$. We proceed to show that
$\mathcal{P}$ is a $\beta $-net in $\mathcal{K}^n$ provided both
$1/\theta $ and $m$ are large enough. In fact, guaranteed by John's
theorem, it is sufficient to prove that $\mathcal{P}$ is a $\beta
$-net in $\overline{\mathcal{K}^n}$.

Assume that $K$ is an $n$-dimensional convex body satisfying
$$B^n\subseteq K\subseteq nB^n.$$
Let ${\bf p}_i$ be the point in $X_{i,m}\cap K$ which is the
furthest to the origin. Then we define
$$P=conv\{ {\bf p}_1, {\bf p}_2, \cdots , {\bf p}_{c(n,\theta )}\}$$
and
$${\bf q}_i=\left\{
\begin{array}{ll}
{\bf p}_i& \mbox{ if ${\bf p}_i={\bf x}_i$}\\
{\bf p}_i+\mbox{${{n-1}\over {mn}}$}{\bf x}_i&
\mbox{otherwise.}\end{array}\right.$$ Clearly the $c(n,\theta )$
caps $({\theta \over n}B^n+{1\over n}{\bf x}_i)\cap \partial (B^n)$
form a covering on $\partial (B^n)$ and therefore we have
$$\left(1-\mbox{$\theta \over n$}\right)B^n\subseteq P\subseteq K.\eqno (4)$$

\begin{figure}[ht]
\centering
\includegraphics[height=4cm,width=7.2cm,angle=0]{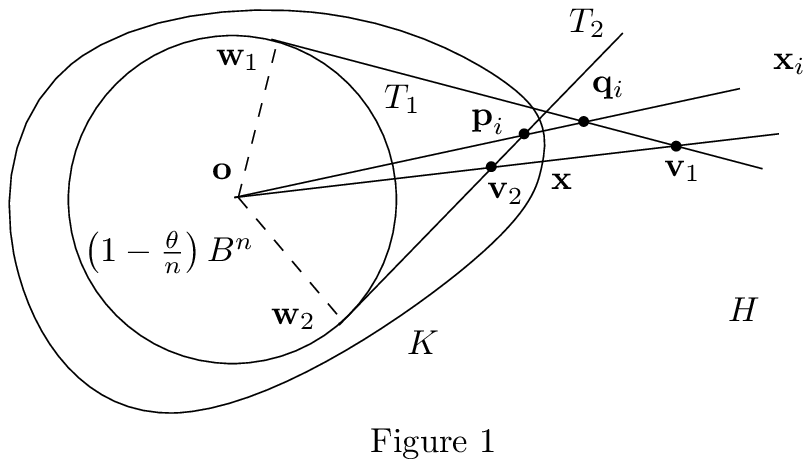}
\end{figure}

Let $C_i$ denote the cone with vertex ${\bf o}$ over the cap
$(\theta B^n+{\bf x}_i)\cap \partial (nB^n)$, let ${\bf x}$ be a
point which is on the boundary of $K$ and belongs to $C_i$, and let
$H$ be a two-dimensional plane passing ${\bf o}$, ${\bf p}_i$ and
${\bf x}$. As shown in Figure 1, $T_1$ is tangent to $(1-{\theta
\over n})B^n\cap H$ at ${\bf w}_1$ and passing ${\bf q}_i$, $T_2$ is
tangent to $(1-{\theta \over n})B^n\cap H$ at ${\bf w}_2$ and
passing ${\bf p}_i$, the straight line determined by ${\bf o}$ and
${\bf x}$ intersects $T_i$ at ${\bf v}_i$. It follows by convexity
that the point ${\bf x}$ is between ${\bf v}_1$ and ${\bf v}_2$. By
elementary geometry, letting $d({\bf u}_1, {\bf u}_2)$ denote the
distance between ${\bf u}_1$ and ${\bf u}_2$, we get
\begin{eqnarray*}
d({\bf v}_1, {\bf v}_2)&\le & d({\bf q}_i, {\bf v}_1)+d({\bf p}_i,
{\bf q}_i)+d({\bf v}_2,{\bf p}_i)\\
&=&(d({\bf w}_1, {\bf v}_1)-d({\bf w}_1,{\bf q}_i))+d({\bf p}_i,{\bf
q}_i)+(d({\bf w}_2, {\bf p}_i)-d({\bf w}_2,{\bf v}_2))\\
&\le & \left(1-\mbox{${\theta \over n}$}\right)\cdot \tan
\left(\arccos \mbox{${{n-\theta }\over
{n^2}}$}+\theta'\right)-\sqrt{n^2-\mbox{$(1-{\theta \over n})^2$}}+\mbox{${n-1}\over m$}\\
&& +\sqrt{\left(n-\mbox{${n-1}\over m$}\right)^2-\mbox{$(1-{\theta
\over n})^2$} }-\left(1-\mbox{${\theta \over n}$}\right)\cdot\tan
\left( \arccos \mbox{${{1-{\theta \over n}}\over {n-{{n-1}\over
m}}}$}-\theta '\right),
\end{eqnarray*}
where $$\theta '=2\arcsin \mbox{$\theta \over {2n}$}.$$ For
convenience, we abbreviate the final complicated function as
$f(n,m,\theta )$. Let ${\bf p}$ denote the point on the boundary of
$P$ and in the direction of ${\bf x}$. Then we have
$${{d({\bf o}, {\bf x})}\over {d({\bf o}, {\bf p})}}\le {{d({\bf o}, {\bf p})+d({\bf v}_1,{\bf v}_2)}\over
{d({\bf o}, {\bf p})}} \le 1+{{f(n,m,\theta )}\over {1-\theta /n}}$$
and therefore
$$K\subseteq \left( 1+\mbox{${{f(n,m,\theta )}\over {1-\theta /n
}}$}\right) P.\eqno (5)$$ In addition, for any fixed $n$, it can be
shown by a routine argument that
$$\lim_{{m\to \infty}\atop {\theta \to 0}}{{f(n,m,\theta )}\over {1-\theta /n }}=0.\eqno (6)$$

As a conclusion of (4), (5) and (6), for any given pair of $n$ and
$\beta $ there is a set $\mathcal{P}$ of $l(n,\beta )$ polytopes
such that for each $K\in \mathcal{K}^n$ there is a $P\in
\mathcal{P}$ satisfying
$$\| K, P\| <\beta .$$
The lemma is proved. \hfill{$\diamondsuit$}

\vspace{0.6cm} Since $$\log (1+x)\le x$$ holds for all $x\ge 0$, it
follows by (4) and (5) that $\mathcal{P}$ will be a $\beta $-net in
$\mathcal{K}^n$ if
$${{f(n,m,\theta )}\over {1-\theta /n}}\le \beta . \eqno (7)$$
By a routine computation one can deduce that, when both $m$ and
$1/\theta $ are sufficiently large and $n$ is fixed,
$$ f(n,m,\theta )\sim 2n \theta +\mbox{${{n-1}\over m}$}<3n\left( \theta +\mbox{${1\over
m}$}\right).$$ Therefore, when $\beta $ is small, (7) can be
guaranteed by taking $\theta ={\beta \over {7n}}$ and $m=\lfloor
{{7n}\over \beta }\rfloor$. In this case, we have $$\theta'=2\arcsin
\mbox{${\beta \over {14n^2}}$}.$$ According to B\"or\"oczky and
Wintsche \cite{bor03} there is a cap covering satisfying
\begin{eqnarray*}
c(n,\theta )&\le  &c\cdot n^{3\over 2}\cos \theta' \sin^{-n}\theta ' \cdot \log (2+n\cos^2\theta')\\
&\le &c\cdot 14^n\cdot n^{2n+3}\cdot \beta^{-n},
\end{eqnarray*}
where $c$ is a suitable constant. Thus, there is a $\beta$-net
$\mathcal{N}=\{ K_1, K_2, \cdots , K_{l(n,\beta )}\}$ in
$\mathcal{K}^n$ satisfying $$l (n, \beta )\le m^{c(n,\theta )}\le
\left\lfloor\mbox{${{7n}\over \beta }$}\right\rfloor^{c\cdot
14^n\cdot n^{2n+3}\cdot \beta^{-n}}.$$

\medskip
\noindent {\bf Remark 4.} Clearly, to estimate the minimal
cardinality of the $\beta $-nets in $\mathcal{K}^n$ and to construct
the corresponding good $\beta $-nets are challenging problems.

\medskip \noindent {\bf The philosophy of our program.}  If
Hadwiger's conjecture is true in $E^n$, since $\gamma_{2^n}(K)$ is
continuous on $\mathcal{K}^n$ and $\mathcal{K}^n$ is compact, then
there is a positive number $c_n< 1$ such that
$$\gamma_{2^n}(K)\le c_n\eqno (8)$$
holds for all $K\in \mathcal{K}^n$. On the other hand, if (8) holds
with certain $c_n<1$ for all convex bodies $K\in \mathcal{K}^n$,
then Hadwiger's conjecture is true in $E^n$. Since $\gamma_{2^n}(K)$
is continuous on $\mathcal{K}^n$, there is a positive number $\beta
$ such that
$$\left| \gamma_{2^n}(K)-\gamma_{2^n}(K')\right| \le \mbox{${1\over 2}$}(1-c_n) \eqno (9)$$
holds whenever $\| K,K'\| \le \beta $. If we can construct a
$\beta$-net $$\mathcal{N}=\left\{ K_1, K_2, \cdots , K_{l(n,\beta
)}\right\}$$ for this particular $\beta $ and can verify (with the
assistance of a computer if necessary) $$\gamma_{2^n}(K_i)\le c_n$$
for all $K_i\in \mathcal{N}$, then by (9) we get
$$\gamma_{2^n}(K)\le \mbox{${1\over 2}$} (1+c_n)$$ for all $K\in
\mathcal{K}^n$ and therefore Hadwiger's conjecture.

\bigskip \noindent {\bf A four-step program for Hadwiger's
conjecture.}

\medskip\noindent
{\bf Step 1.} In the considered dimension, for example $n=3$, study
the values of $\gamma_{2^n}(K)$ for some particular convex bodies
$K$ and therefore choose a possible candidate constant $c_n$ for
(8).

\medskip\noindent
{\bf Step 2.} Choose a suitable positive number $\beta $ to
guarantee (9), based on a close study on the function
$\gamma_{2^n}(K)$.

\medskip\noindent
{\bf Step 3.} Construct a suitable $\beta$-net $\mathcal{N}$ based
on the fundamental lemma.

\medskip\noindent
{\bf Step 4.} Verify that
$$\gamma_{2^n}(K_i)\le c_n$$
holds for all $K_i\in \mathcal{N}$ (with the assistance of a
computer if necessary).

\medskip\noindent
{\bf Remark 5.} In principle, the conjecture can be proved in $E^n$
by our program if it is true in this particular dimension and if the
computing facility is efficient enough. Clearly the set
$\mathcal{P}$ can be much reduced in cardinality.

\bigskip
Let $d(X)$ denote the {\it diameter} of a bounded set $X$ in
$E^n$. In other words,
$$d(X)=\sup \{ d({\bf x}_1, {\bf x}_2):\ {\bf x}_i\in X\},$$
where $d({\bf x}_1, {\bf x}_2)$ denotes the Euclidean distance
between ${\bf x}_1$ and ${\bf x}_2$. In 1933, Borsuk \cite{bors33}
made a conjecture that {\it each bounded $n$-dimensional set $X$ can
be partitioned into $n+1$ parts $X_1$, $X_2$, $\cdots $, $X_{n+1}$
such that}
$$d(X_i)<d(X), \qquad i=1, 2, \cdots , n+1.$$
This conjecture has attracted much attention. The two-dimensional
case was proved by Borsuk himself, the three-dimensional case was
first proved by J. Perkal \cite{perk47} in 1947, and the case of
sufficiently high dimensions was disproved by J. Kahn and G. Kalai
\cite{kahn93} in 1993. For a survey on this conjecture we refer to
Zong \cite{zong08}. Up to now, it is still open for $4\le n\le 297$.

It is known in convexity that for each bounded set $X$ there is a
convex body $\widehat{X}$ of constant width satisfying both
$X\subseteq \widehat{X}$ and
$$d(X)=d(\widehat{X}).$$
Therefore, to study Borsuk's conjecture it is sufficient to consider
all convex bodies of constant width $1$.

Let $r(K)$ and $R(K)$ denote the radii of the maximal {\it insphere}
and the minimal {\it circumsphere} of an $n$-dimensional convex body
$K$ of constant width $1$, respectively. It is known in convexity
(see \cite{chak83}) that
$$1-\sqrt{{n\over {2n+2}}}\le r(K)\le R(K)\le \sqrt{{n\over
{2n+2}}}.$$ For convenience, we take
$$\mu_n=\sqrt{{n\over {2n+2}}} \left/\left( 1-\sqrt{{n\over {2n+2}}}
\right)\right.={{\sqrt{2n^2+2n}+n}\over {n+2}}$$ and define
$\widehat{\mathcal{K}^n}$ to be the set of all $n$-dimensional
convex bodies $K$ satisfying $B^n\subseteq K\subseteq \mu_n B^n$
associated with the {\it Hausdorff metric} $\| \cdot \|^*$, where
$$\| K_1, K_2\|^*=\min\{ r:\ K_1\subseteq K_2+rB^n, \ K_2\subseteq K_1+rB^n\}.$$
Then, we define $\varphi_m (K)$ to be the minimal number $r$ such
that $K$ can be partitioned into $m$ parts $K_1$, $K_2$, $\cdots $,
$K_m$ such that
$$d(K_i)\le r\cdot d(K)$$
holds for all $i=1,$ $2,$ $\cdots $, $m$.

It is easy to see that $\widehat{\mathcal{K}^n}$ is compact and, for
any positive integer $m$,  $\varphi_m (K)$ is continuous on
$\widehat{\mathcal{K}^n}$. To prove Borsuk's conjecture, it is
sufficient to show
$$\varphi_{n+1} (K)\le c_n<1,\qquad K\in \widehat{\mathcal{K}^n}.$$
Therefore, in a given dimension, for instance $n=4$, Borsuk's
conjecture should be approachable by a four-step program similar to
that for Hadwiger's conjecture.

\medskip\noindent
{\bf Remark 6.} It was shown (see Gr\"unbaum \cite{grun63}) that
$$\varphi_3(K)\le {{\sqrt{3}}\over 2}$$
holds for all two-dimensional $K$, and
$$\varphi_4(K)\le 0.9887$$
holds for all three-dimensional $K$.

\section*{4.  The covering functions on $\mathcal{K}^3$}

\medskip
In this section, among other things, we will prove Theorem 1 and
Theorem 2. As a consequence, we give some insight to a reasonable
estimate for the constant $c_3$ defined in the previous section.
First, let us introduce two lemmas.

\medskip\noindent
{\bf Lemma 1 (Besicovitch \cite{besi48}).} {\it Each two-dimensional
convex domain has an inscribed affine regular hexagon. }

\medskip\noindent
{\bf Remark 7.} Affine regular hexagons are the imagines of a
regular hexagon under non-singular linear transformations.

\medskip\noindent
{\bf Lemma 2.} {\it Let $K$ be a two-dimensional convex domain, let
$\lambda $ be a real number satisfying $0<\lambda <1$, and let $\{
{\bf x}_1, {\bf x}_2, {\bf x}_3\}$ be an ordered triple on the
boundary of $K$. If $\{ {\bf x}_1, {\bf x}_2, {\bf x}_3\}\subset
\lambda K +{\bf y}$, then the whole curve from ${\bf x}_1$ to ${\bf
x}_3$ passing ${\bf x}_2$ belongs to $\lambda K+{\bf y}$.}

\medskip\noindent
{\bf Proof.} For convenience, we assume that ${\bf o}\in int(K)$ and
${\bf y}=(0, a)$. It is well known in convexity (see Eggleston
\cite{eggl58}) that the set of {\it regular convex domains} (each
tangent touches $K$ at exactly one point and there is one and only
one tangent at each boundary point) is dense in $\mathcal{K}^2$.
Therefore, without loss of generality, we assume that $K$ is
regular.

\begin{figure}[ht]
\centering
\includegraphics[height=5cm,width=7cm,angle=0]{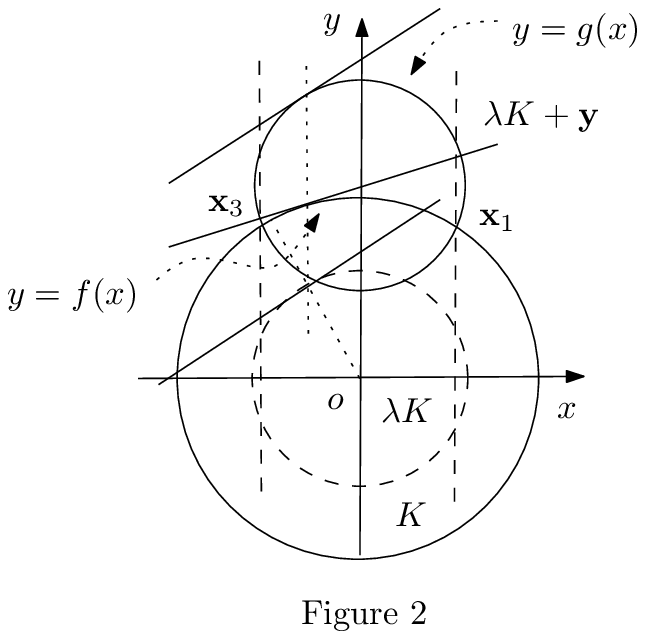}
\end{figure}

Let ${\bf x}_1=(x_1, y_1)$ and ${\bf x}_3=(x_3, y_3)$ denote the
points of $\partial (K)\cap (\lambda K+{\bf y})$ with maximal and
minimal $x$-coordinates, respectively. Let $y=f(x)$ denote the curve
of $\partial (K)$ from ${\bf x}_3$ to ${\bf x}_1$, and let $y=g(x)$
denote the above part of $\partial (\lambda K)+{\bf y}$ in the strip
of $x_3\le x\le x_1$. By convexity, as shown in Figure 2, we have
$$g(x_3)\ge f(x_3),$$
$$g(x_1)\ge f(x_1),$$
$$g'(x)=f'(\mbox{${1\over \lambda}$}x )\ge f'(x)$$
for $ x_3\le x\le 0$, and
$$g'(x)=f'(\mbox{${1\over \lambda}$}x)\le f'(x)$$
for $0\le x\le x_1.$ Thus, we get
$$g(x)-f(x)=g(x_3)-f(x_3)+\int_{x_3}^x(g'(t)-f'(t))dt\ge 0$$
when $x_3\le x\le 0$, and
$$g(x)-f(x)=g(x_1)-f(x_1)+\int_0^x(f'(t)-g'(t))dt\ge 0$$
when $0\le x\le x_1$. Therefore, by convexity, the whole curve
$y=f(x)$ belongs to $\lambda K+{\bf y}$. The lemma is proved.
\hfill{$\diamondsuit $}

\medskip\noindent
{\bf Corollary 1.} {\it Let $K$ be an $n$-dimensional convex body,
$\lambda $ be a real number satisfying $0<\lambda <1$, $R$ be a
closed region on $\partial (K)$ with boundary $\Gamma $ and a
relatively interior point ${\bf p}$. If $$\Gamma \cup \{ {\bf p}\}
\subset \lambda K+{\bf y}$$ holds for some point ${\bf y}$, then we
have} $$R\subset \lambda K+{\bf y}.$$

\medskip\noindent
{\bf Proof of Theorem 1.} Let $K$ be a three-dimensional cone over a
convex domain $D$. By Lemma 1, there is an affine regular hexagon
$H$ inscribed in $D$. Without loss of generality, we assume that
${\bf v}=(0,0,1)$ is the vertex of $K$, $H$ is perpendicular to
${\bf v}$ and centered at the origin ${\bf o}=(0,0,0)$.

\begin{figure}[ht]
\centering
\includegraphics[height=5cm,width=8cm,angle=0]{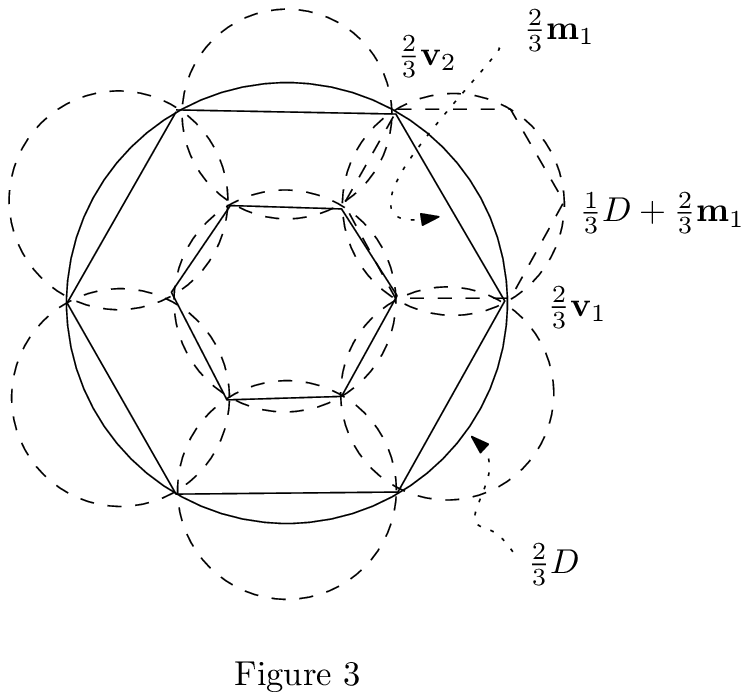}
\end{figure}

\begin{figure}[ht]
\centering
\includegraphics[height=4.5cm,width=7.3cm,angle=0]{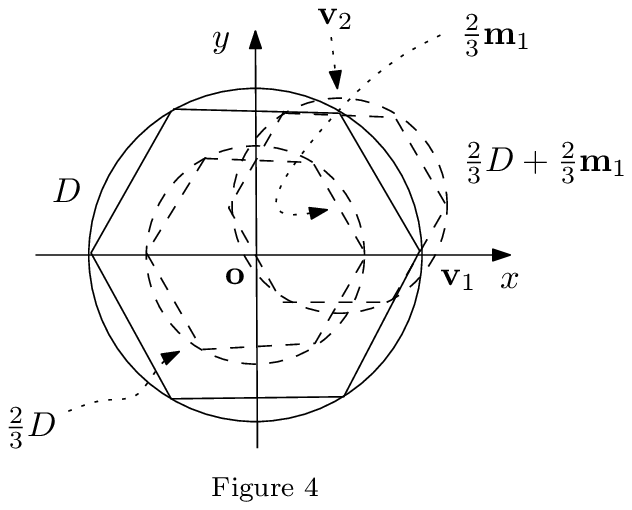}
\end{figure}

Let ${\bf v}_1,$ ${\bf v}_2$, $\cdots $, ${\bf v}_6$ be the six
vertices of $H$ and let ${\bf m}_1$, ${\bf m}_2$, $\cdots $, ${\bf
m}_6$ denote the midpoints of ${\bf v}_1{\bf v}_2$, ${\bf v}_2{\bf
v}_3$, $\cdots$, ${\bf v}_6{\bf v}_1$, respectively. By elementary
argument, as shown in Figure 3, we have
$$\left\{ \mbox{${2\over 3}$}{\bf v}_1, \mbox{${2\over 3}$}{\bf v}_2\right\} \subset
\mbox{${1\over 3}$}D+\mbox{${2\over 3}$}{\bf m}_1.$$ Thus, by Lemma
2 we get
$$\mbox{${2\over 3}$}D+\mbox{${1\over 3}$}{\bf v}\subseteq
\bigcup_{i=0}^6\left(\mbox{${1\over 3}$}D+\mbox{${2\over 3}$}{\bf
m}_i+\mbox{${1\over 3}$}{\bf v}\right), \eqno (10)$$ where ${\bf
m}_0=(0,0,0)$. Similarly, as shown in Figure 4, we have
$$\{ {\bf v}_1, {\bf v}_2\}\subset \mbox{${2\over 3}$}D+\mbox{${2\over 3}$}{\bf m}_1$$
and therefore
$$D\subseteq \bigcup_{i=0}^6\left(\mbox{${2\over 3}$}D+\mbox{${2\over 3}$}{\bf
m}_i\right).\eqno (11)$$ On the other hand, we have $$\mbox{${1\over
3}$}D+\mbox{${2\over 3}$}{\bf m}_i+\mbox{${1\over 3}$}{\bf v}
\subset \mbox{${2\over 3}$}K+\mbox{${2\over 3}$}{\bf m}_i\eqno
(12)$$ and $$\mbox{${2\over 3}$}D+\mbox{${2\over 3}$}{\bf
m}_i\subset \mbox{${2\over 3}$}K+\mbox{${2\over 3}$}{\bf m}_i.\eqno
(13)$$ Therefore, by (10), (11), (12), (13) and convexity we get
$$K\subseteq \bigcup_{i=0}^7\left( \mbox{${2\over 3}$}K+\mbox{${2\over 3}$}{\bf
m}_i\right),$$ where ${\bf m}_7={1\over 2}{\bf v}.$ Theorem 1 is
proved. \hfill{$\diamondsuit $}

\bigskip\noindent {\bf Proof of Theorem 2.}

\medskip\noindent
{\bf Case 1.} $1\le p\le 2$. In this case, we take
$$\Gamma=\left\{ {\bf x}=(x_1, x_2, x_3):\ x_1=\left(\mbox{${1\over
3}$}\right)^{1\over p},\ {\bf x}\in \partial (K_p)\right\},$$
$\lambda =\sqrt{2\over 3}$, ${\bf y}=(\left({1\over
3}\right)^{1\over p}, 0, 0)$, and let $R$ denote the part of
$\partial (K_p)$ bounded by $\Gamma $ and containing $(1,0,0)$.

For any point ${\bf x}\in \Gamma $, we have
$$\left(\left(\mbox{$1\over 3$}\right)^{1\over
p}\right)^p+|x_2|^p+|x_3|^p=1,$$
$$|x_2|^p+|x_3|^p=\mbox{${2\over 3}$}\le \left(\sqrt{\mbox{$2\over
3$}}\right)^p$$ and therefore
$$\Gamma \subset \lambda K_p+{\bf y}.$$
On the other hand, it can be verified that
$$(1,0,0)\in \lambda K_p+{\bf y}.$$
By Corollary 1 we have $$R\subset \lambda K_p+{\bf y}.$$ Therefore,
in this case $K_p$ can be covered by the union of $\lambda K_p\pm
\left( ({1\over 3})^{1\over p}, 0,0\right)$, $\lambda K_p\pm \left(
0, ({1\over 3})^{1\over p},0\right)$ and $\lambda K_p\pm \left( 0,
0, ({1\over 3})^{1\over p}\right)$ and thus
$$\gamma_8(K_p)\le \gamma_6(K_p)\le \sqrt{\mbox{$2\over 3$}}.$$

\noindent {\bf Case 2.} $2\le p\le \infty $. In this case we define
$$\Gamma_i=\left\{{\bf x}=(x_1, x_2,x_3):\ x_i=0,\ x_j\ge 0, j\not= i,\
{\bf x}\in \partial (K_p)\right\},$$
$$\Gamma =\Gamma_1\cup \Gamma_2\cup\Gamma_3,$$
$\lambda =\sqrt{2\over 3}$, ${\bf y}=({1\over 3}, {1\over 3},
{1\over 3})$, and let $R$ denote the part of $\partial (K_p)$
bounded by $\Gamma $ and containing the point $(({1\over 3})^{1\over
p}, ({1\over 3})^{1\over p},({1\over 3})^{1\over p})$.

Let $J$ denote the intersection of $K_p$ with the plane $x_1=0$, and
let $J'$ denote the intersection of $\lambda K_p+{\bf y}$ with the
plane. It is easy to see that $J'$ is homothetic to $J$. By a
routine computation, for all $2\le p\le +\infty$, it can be shown
that
$$\left( \mbox{$2\over 3$}\right)^p +2\left( \mbox{$1\over
3$}\right)^p\le \left( \mbox{$2\over 3$}\right)^{p\over 2}.$$ Thus,
both $(0,1,0)$ and $(0,0,1)$ belong to $J'$. Consequently, we also
have $$\left( 0, (\mbox{$1\over 2$})^{1\over p}, (\mbox{$1\over
2$})^{1\over p}\right)\in J'.$$ By Lemma 2, we get
$$\Gamma_1\subset J'\subset \lambda K_p+{\bf y}$$
and therefore
$$\Gamma \subset \lambda K_p+{\bf y}.$$
On the other hand, it can be verified that $${\bf o}\in \lambda
K_p+{\bf y},$$ $$3\left( \left( \mbox{$1\over 3$}\right)^{1\over
p}-\mbox{$1\over 3$}\right)^p\le \left( \mbox{$2\over
3$}\right)^{p\over 2}$$ and therefore
$$\left(\left(\mbox{$1\over 3$}\right)^{1\over p}, \left(\mbox{$1\over 3$}\right)^{1\over p},
\left(\mbox{$1\over 3$}\right)^{1\over p}\right)\in \lambda K_p+{\bf
y}.$$ By Corollary 1 we get
$$R\subset \lambda K_p+{\bf y}.$$
Thus, in this case $K_p$ can be covered by the union of the eight
translates $\sqrt{2\over 3}K_p+\left( \pm \mbox{$1\over 3$}, \pm
\mbox{$1\over 3$}, \pm \mbox{$1\over 3$}\right)$ and hence
$$\gamma_{8}(K_p)\le \sqrt{\mbox{$2\over 3$}}.$$

As a conclusion of the two cases, Theorem 2 is proved.
\hfill{$\diamondsuit $}

\bigskip\noindent {\bf Remark 8.} It was shown by Sch\"utte
\cite{schu55} that
$$\gamma_8(K_2)\le \sin 48^\circ 9'=0.744894\cdots <\sqrt{\mbox{$2\over
3$}}.$$ Thus, it follows from the proof of Theorem 2 that
$\sqrt{2\over 3}$ is not the optimal upper bound for
$\gamma_8(K_p).$ However, perhaps one can take $c_3=\sqrt{2\over
3}$.

\bigskip\noindent {\bf Remark 9.} Let $T$ denote a regular
tetrahedron. The values of $\gamma_m(K)$ for some small $m$ and some
particular $K$ are listed in Table 3. The values of $\gamma_4(K_2)$
and $\gamma_6(K_2)$ were discovered by L. Fejes T\'oth \cite{feje53}
and the values of $\gamma_5(K_2)$ and $\gamma_7(K_2)$ were
determined by Sch\"utte \cite{schu55}.

\begin{table}
\renewcommand\arraystretch{1.5}
\noindent\[
\begin{array}{|c|c|c|c|c|c|}
\hline
m&4&5&6&7&8\\
\hline
\gamma_m (T)&{3\over 4}&{9\over 13}&?&?&?\\
\hline
\gamma_m (K_1)&1&1&{2\over 3}&{2\over 3}&{2\over 3}\\
\hline
\gamma_m (K_2)&0.9428\cdots &0.8944\cdots &0.8164\cdots &0.7775\cdots &?\\
\hline
\end{array}
\]
\centerline{Table 3}
\end{table}

\bigskip\noindent {\bf Remark 10.} If Hadwiger's conjecture is true
for all dimensions, then we have
$$\Gamma (n,2^n)<1$$
for all $n$. Nevertheless, it seems that
$$\lim_{n\to\infty }\Gamma (n,2^n)=1.$$

\bigskip\medskip\noindent
{\large\bf Acknowledgement.} I am very grateful to Professor Peter
M. Gruber, Professor Martin Henk and the referees for their helpful
comments and remarks.

\bibliographystyle{amsplain}

\end{document}